\documentclass[12pt]{amsart}
\usepackage{graphicx}
\usepackage{amssymb}
\usepackage{epstopdf}
 \usepackage{latexsym}
\usepackage{amsmath}
\usepackage{amsfonts} 
\usepackage{amsthm}

\newtheorem{Proposition}{Proposition}

\newtheorem{Lemma}[Proposition]{Lemma}

\newtheorem{Corollary}[Proposition]{Corollary}

\def\XXint#1#2#3{{\setbox0=\hbox{$#1{#2#3}{\int}$}
\vcenter{\hbox{$#2#3$}}\kern-.5\wd0}}

\def\e{\epsilon}

    \def\sqr#1#2{{\vcenter{\vbox{\hrule height .#2pt
                             \hbox{\vrule width .#2pt height#1pt \kern#1pt
                                   \vrule width .#2pt}
                             \hrule height .#2pt}}}}

     \def\CC{\mathbb{C}}

    \def\e{\epsilon}

\def\be{\begin{equation}}
\def\ee{\end{equation}}

\def\al{\alpha}

\begin{document}

\title[First return map ]{The return map for a planar vector field with nilpotent linear part: a direct and explicit derivation}
\author{Rodica D. Costin}

\maketitle

\begin{abstract}

Using a direct approach the return map near a focus of a planar vector field with nilpotent linear part is found as a convergent power series which is a perturbation of the identity and whose terms can be calculated iteratively. The first nontrivial coefficient is the value of an Abelian integral, and the following ones are explicitly given as iterated integrals.

\end{abstract}


\section{Introduction}

The study of planar vector fields has been the subject of intense research, particularly in connection to Hilbert's 16th Problem. Significant progress has been made in the geometric theory of these fields, as well as in bifurcation theory, normal forms, foliations,  and the study of Abelian integrals \cite{Dum_Book},\,\cite{Rouss_3}.

The Poincar\'e first return maps have been studied in view of their relevance for establishing the existence of closed orbits, and also due to their large number of applications (see for example \cite{Guckenheimer} and references therein), and also in connection to o-minimality \cite{Miller}.

A fundamental result concerns the asymptotic form of return maps states that if the singular points of a $C^\infty$ vector field are algebraically isolated, there exists a semitransversal arc such that the return map admits an asymptotic expansion is positive powers of $x$ and logs (with the first term linear), or has its principal part a finite composition of powers and exponentials \cite{Ilyashenko},\,\cite{Moussu}.

In the case when the linear part of the vector field has non-zero eigenvalues there are important results containing the return map \cite{Brudnyi},\,\cite{Dumortier},\,\cite{Rouss_1},\,\cite{Rouss_2}, \cite{Francoise_Lienard_2},\,\cite{Francoise_Lienard_1},\,\cite{ret1}. Results are also available for perturbations of Hamiltonians \cite{Francoise_Iterated_int}, \cite{Francoise_Suc_der} and for perturbations of integrable systems \cite{Gavrilov}.
On the other hand, there are few results available in the general setting \cite{Berezovskaya},\,\cite{Medvedeva_Lt},\,\cite{Medvedeva_St}.

The present paper studies an example of a field with nilpotent linear part, near a focus. The main goal is to establish techniques that allow to deduce the return map as a suitable series which can be calculated algorithmically and can be used in numerical calculations.

 \section{Main result}

The paper studies the return map for the system
\be\label{sysXY}
\dot{X}=-Y,\ \ \ \ \dot{Y}=X^3-Y^3
\ee
which has a nilpotent linear part (both eigenvalues are zero). This is one of the simplest examples of systems in this class \cite{Perko}, and for which there are (to the author's knowledge) no methods available to generate the return map.

The main result is the following:

\begin{Proposition}\label{MainProp}

Let  $\epsilon$ with $0<\epsilon<\epsilon_0$ small enough.

The solution of (\ref{sysXY}) satisfying $X(0)=\epsilon,\,Y(0)=0$ first returns to the positive $X$-axis at the value $\tilde{\tilde{X}}$ satisfying
$$\tilde{\tilde{X}}=\epsilon+\sum_{n=1}^\infty X_n\epsilon^{3n+1}$$
which is a convergent series.

The coefficients $X_n$ can be calculated iteratively.
In particular,
$$
X_1=-2^{3/2}\,c_{{1}},\ \ X_2=16\,{c_{{1}}}^{2},\ \ X_3=-2^{3/2}\, \left( 20\,c_{{2}}c_{{1}}+8\,c_{{3}}+49\,{c_{{1}}}^{3}\right)$$
where $c_n=v_n(1)$ with $v_n$ given by
\be\label{v12xi}
v_1(\xi)=\xi^{-2}\int_0^\xi t^4p(t)^{3/2}\, dt,\ \ \ v_2(\xi)=-\frac{3}{2}\,\xi^{-2}\int_0^\xi t^4p(t)^{1/2}v_1(t)\, dt
\ee
\be\label{v3xi}
v_3(\xi)=\xi^{-2}\int_0^\xi t^4\, \left[ -\frac {3}{2}\,
{p(t)}^{1/2}v_{{2}}(t)+\frac {3}{8}\,{\frac {{v_{{1}}}^{2}(t)}{{p}(t)^{1/2}}} \right]\, dt
\ee
where
\be\label{defp}
p(t)=4 -6\,{t}^{2}+4\,{t}^{4}-{t}^{6}
\ee

\end{Proposition}

 \section{Proof of Proposition\,\ref{MainProp}}
 
The proof of Proposition\,\ref{MainProp} also provides an algorithm for calculating iteratively the coefficients $X_n$.

\subsection{Normalization}\label{normalization}

It is convenient to normalize the variables $X,Y,t$ so that the constant $\e$ appears as a small parameter in the equation: with 
$$X=\e x,\ \ Y=2^{-1/2}\e^2 y,\ \tau=2^{-3/2}\e t$$
the system (\ref{sysXY}) becomes
\be\label{syse}
\frac{dx}{d\tau}=-2y,\ \frac{dy}{d\tau}=4x^3-\alpha y^3
\ee
where $\al$ is the small parameter:
\be\label{not_alpha}
\alpha=2^{1/2}\e^3
\ee

While the initial condition $X(0)=\e$ becomes $x(0)=1$, it is useful to study solutions of (\ref{syse}) with the more general initial condition $x(0)=\eta$ with $\eta$ in a neighborhood of $1$.

{\bf{Remark.}} The system (\ref{syse}) has the form $dH+\alpha\omega=0$ with $H=y^2+x^4$, $\omega=y^3\, dx$ and $\alpha$ a small parameter. The recent result \cite{Francoise_Iterated_int} gives a formalism for finding the return map for this type of systems; see also \cite{Gavrilov}. The present construction is concrete, explicit, and suitable for numerical calculations.

\subsection{General behavior of solutions of (\ref{syse}).} 

Let $\al>0$. 

Note the following Lyapunov function for (\ref{syse}):
\be\label{Lyapunov}
L(x,y)=y^2+x^4,\ \  {\mbox{with}}\ \ \frac{d}{d\tau}L\left(x(\tau),y(\tau)\right)=-2\al \,y(\tau)^4\leq 0
\ee

Since the set $\{y=0\}$ contains no trajectories besides the origin (which is the only equilibrium point of (\ref{syse})), then the origin is asymptotically stable by the Krasovskii-LaSalle principle.

Consider the solution of (\ref{syse}) with the initial condition $x(0)=\eta,\, y(0)=0$ for some $\eta\in[1/2,3/2]$. 

Since $y'(0)>0$ and $x'(0)=0,\,x''(0)<0$ then $y$ increases and $x$ decreases for small $\tau>0$. This monotony must change due to (\ref{Lyapunov}), and this can happen only at some point where $y=0$ or where $4x^3=\al y^3$, whichever comes first. Since $y$ increases, then the first occurrence is a point where $4x^3=\al y^3$. At this point $x'<0$ so $x$ continues to decrease, while $y''=-24x^2y<0$ so $y$ has a maximum, and will continue by decreasing. Again, the monotony must change due to (\ref{Lyapunov}) and the path cannot cross again the line $4^{1/3}x=\al^{1/3}y$ before the monotony of $x$ changes, therefore the the next change of monotony happens for $y=0$, a point where, therefore $x<0$; denote this value of $x$ by $-\tilde{\eta}$.

Denote by $x_0(\tau),y_0(\tau)$ the solution for $\al=0$: $x_0'=-2y,\, y_0'=4x_0^3$ and $x_0(0)=\eta,\,y_0(0)=0$. Therefore $x_0^4+y_0^2=\eta^4$. We have $x(\tau)=x_0(\tau)+O(\al),\, y(\tau)=y_0(\tau)+O(\al)$ therefore $\tilde{\eta}=\eta+O(\al)$.

Similar arguments show that the solution $x(\tau),y(\tau)$ of (\ref{syse}) continues to turn around the origin, crossing again the positive $x$-axis at a point $\tilde{\tilde{\eta}}=\eta+O(\al)$ (and, of course, $\tilde{\tilde{\eta}}<\eta$ by (\ref{Lyapunov})).

Solutions $(x(\tau),y(\tau))$ of (\ref{syse}) provide smooth parametrizations for solutions $y(x)$ of
\be\label{eqy}
\frac{d}{dx}\left(y^2\right)=-4x^3+\al y^3
\ee

Note that following the path $(x(\tau),y(\tau))$ one full rotation around the origin corresponds to considering a positive solution of (\ref{eqy}), followed by a negative one.

\subsection{Positive solutions of (\ref{eqy}) for $x>0$.}\label{QI}

Lemma\,\ref{L1} shows that there exists a unique solution $y\geq 0$ of (\ref{eqy}) so that $y(\eta)=0$, that this solution is defined for $x\in[0,\eta]$ and establishes an iterative procedure for calculating this solution.

Substituting $y=u^{1/2}$ in equation (\ref{eqy}) we obtain
\be\label{equ}
\frac{du}{dx}=-4x^3+\al u^{3/2}
\ee

\begin{Lemma}\label{L1}

There exists $\delta_0>0$ independent of $\eta$ and $\al$ so that so that the following holds.

Let $\eta>0$. For any $\al$ with $|\al |<\al_0\leq \delta_0\eta^{-3}/2$ equation (\ref{equ}) with the condition $u(\eta)=0$ has a unique solution $u=u(x;\al,\eta)$ for $x\in [0,\eta]$. We have $u(x;\al,\eta)>0$ for $x\in[0,\eta)$, $\al>0$ and $u(x;\al,\eta)$ is analytic in $\al$ and $\eta$ for $|\al |<\al_0$ and $\eta>0$.

\end{Lemma}

{\bf{Remark.}} We will use the results of Lemma\,\ref{L1} only for  $\eta$ such that $1/2\leq\eta\leq3/2$. In this case (by lowering $\al_0$) we can take $\al_0$ of Lemma\,\ref{L1} independent of $\eta$ by taking
\be\label{estimal0}
\al_0\leq \delta_0/2\min \, \eta^{-3}=4^3/{27^3}\delta_0/2\equiv c_0\delta_0
\ee

{\em{Proof of Lemma\,\ref{L1}.}}

Local analysis shows that solutions of (\ref{equ}) satisfying $u(\eta)=0$ have an expansion in integer and half-integer powers of $\eta-x$ and we have
\be\label{expuv}
u(x)=\eta^4-x^4-\frac{16}{5}\,\alpha\,\eta^{9/2}(\eta-x)^{5/2}\left(1+o(1)\right)\ \ \ (x\to\eta-)
\ee
which inspires the following substitutions.

Denote
$$\xi=\sqrt{1-\frac{x}{\eta}}$$
(with the usual branch of the square root for $x/\eta<1$)
and let
$$u \left( x \right) ={\eta}^{4}-{x}^{4}-{\eta}^{3} \left( \eta-x \right) v
 \left( \xi \right) $$
 
 Note that if $u(\eta)=0$ then necessarily $v(0)=0$ by (\ref{expuv}).

Equation (\ref{equ}) becomes
\be\label{eqv}
\xi\,{\frac {d}{d\xi}}v \left( \xi \right) +2\,v \left( \xi \right) =
\delta \,{\xi}^{3}\,\left[\,p(\xi)-v \left( \xi \right)\,  \right] ^{3/2}
\ee
where
\be\label{defdelta}
\delta=2\,{\eta}^{3}\alpha
\ee
and the polynomial $p$ is given by (\ref{defp}).
(Note that $p(\xi)\in[1,4]$ for $\xi\in[0,1]$.)

Lemma\,\ref{L1} follows if we show the following:

\begin{Lemma}\label{Lserv}

These exists $\delta_0>0$ so that for any $\delta$ with $|\delta |<\delta_0$ equation (\ref{eqv}) has a unique solution $v=v(\xi;\delta)$ for $\xi\in[0,1]$ so that $v(0)=0$. 

Moreover, $v(\xi;\cdot)$ is analytic for $\delta\in\CC$ with $|\delta |<\delta_0$ and the terms of its power series 
\be\label{servxi}
 v(\xi;\delta)=\sum_{n\geq 1}\delta^nv_n(\xi)
 \ee
 can be calculated recursively; in particular, the first terms are (\ref{v12xi}),\,(\ref{v3xi}).

\end{Lemma}

To prove Lemma\,\ref{Lserv} multiply (\ref{eqv}) by $\xi$ and integrate; we obtain that $v$ is a fixed point ($v=\mathcal{J}[v]$) for the operator
\be\label{defJ}
 \mathcal{J}v\,\left( \xi \right) =\delta \xi^{-2}\int_0^\xi{t}^{4}\left[\,p(t)-v \left( t \right)\,  \right] ^{3/2}dt=\delta \,\xi^{3}\int_0^1{s}^{4}\left[\,p(\xi s)-v \left( \xi s\right)\,  \right] ^{3/2} ds
 \ee
Let $\mathcal{B}$ be the Banach space of functions $f(\xi;\delta)$ continuous for $\xi\in [0,1]$ and analytic on the (complex) disk $|\delta|<\delta_0$, continuous on $|\delta|\leq\delta_0$, with the norm 
$$\|f\|=\sup_{\xi\in[0,1]}\,\sup_{|\delta|\leq\delta_0}\, |f(\xi;\delta)|$$

Let $m$ be a number with $0<m<1$. Let $\delta_0>0$ be small enough, so that $\delta_0<m/\sqrt{5}$ and $\delta_0<10/3/\sqrt{5}$.

Let $\mathcal{B}_m$ be the ball $\mathcal{B}_m=\{ f\in \mathcal{B}_m\, ;\, \|f\|\leq m\}$. 

We have $\mathcal{J}:\mathcal{B}_m\to\mathcal{B}_m$. Indeed, note that for $f\in \mathcal{B}_m$ we have
$$|p(t)-f(t;\delta)|\geq |p(t)|-| f(t;\delta)|\geq 1-m>0$$
therefore, since $f$ is analytic in $\delta$, then so is $\left[\,p(t)-f \left( t;\delta \right)\,  \right] ^{3/2}$, and therefore so is $\mathcal{J}f$. Also, if $f\in \mathcal{B}_m$ then also $\mathcal{J}f\in \mathcal{B}_m$ because
$$\big| \mathcal{J}f\,(\xi;\delta) \big|\leq  |\delta| \int_0^1 s^4 \left[ |p( \xi s)|+|f( \xi s;\delta)|\right]^{3/2}ds\leq \frac{\delta_0(4+m)^{3/2}}{5}<\delta_0\sqrt{5}<m$$

Moreover, the operator $\mathcal{J}$ is a contraction on $\mathcal{B}_m$. Indeed, using the estimate
$$\big| (p-f_1)^{3/2}-(p-f_2)^{3/2}\big|\leq |f_1-f_2|\, \frac{3}{2}\,\sup_{|f|\leq m}|p-f|^{1/2}\leq |f_1-f_2|
\frac{3(4+m)^{1/2}}{2}$$
we obtain
$$\big| \mathcal{J}f_1-\mathcal{J}f_2 \big|\leq c\|f_1-f_2\|\ \ \ {\mbox{with}}\ c=\delta_0\frac{3\sqrt{5}}{10}<1$$

Therefore the operator $\mathcal{J}$ has a unique fixed point, which is the solution $v(\xi;\delta)$. 

To obtain the power series (\ref{servxi}) substitute an expansion $ v(\xi;\delta)=v_0(\xi)+\sum_{n\geq 1}\delta^nv_n(\xi)$ in (\ref{eqv}). It follows that $\xi v_0'+2 v_0=0$ with $v_0(0)=0$, therefore $v_0(\xi)\equiv 0$.

Substitution of (\ref{servxi}) in $(p-v)^{3/2}$ followed by expansion in power series in $\delta$ give
$$(p-v)^{3/2}=p^{3/2}\left(1-\sum_{n\geq 1}\delta^n\frac{v_n}{p}\right)^{3/2}\equiv p^{3/2}\left(1+\sum_{n\geq 1}\delta^nR_n\right)$$
where $R_n=R_n[v_1,\ldots,v_n,p]$. In particular,
$$R_1=-\frac {3}{2}\,\frac{v_{{1}}}{p},\ \ \ R_2= \left( -\frac {3}{2}\,\frac{v_{{2}}}{p}+\frac {3}{8}\,{\frac {{v_{{1}}}^{2}}{{p^2}}} \right) $$

From (\ref{eqv}) we obtain the recursive system
$$\xi\,{\frac {d}{d\xi}}v_n +2\,v_n  ={\xi}^{3}\, {p}^{{3/2}}R_{n-1}$$
(for $n\geq 1$ and with $R_0=1$), with the only solution with $v_n(0)=0$ given recursively by
$$v_n(\xi)=\xi^{-2}\int_0^\xi\, t^4\,p(t)^{3/2}\,R_{n-1}(t)\, dt$$

In particular, we have (\ref{v12xi}), (\ref{v3xi}).
\qed

\

The following gathers the conclusions of the present section.
\begin{Corollary}\label{Coro1}

There exists $\al_0>0$ so that for any $\eta\in[1/2,3/2]$ and $\al$ with $|\al |<\al_0$ equation (\ref{eqy}) has a unique solution $y(x)$ on $[0,\eta]$ satisfying $y(\eta)=0$ and $y>0$ on $[0,\eta)$ for $\al>0$. 

Moreover, this solution has the form
\be\label{formphi}
y=\phi(x;\al,\eta)=\left[ {\eta}^{4}-{x}^{4}-{\eta}^{3} \left( \eta-x \right) v
 \left( \sqrt{1-\frac{x}{\eta}} ;2\,{\eta}^{3}\alpha \right) \right] ^{1/2}
 \ee
with $v=v(\xi,\delta)$ a solution of (\ref{eqv}). The map $\al\mapsto v(\xi;2\,{\eta}^{3}\alpha )$ is analytic for $|\al|<\al_0$.

\end{Corollary}

\qed

\newpage

\subsection{Solutions of (\ref{eqy}) in other quadrants and matching}\label{QA}

\subsubsection{Solutions in other quadrants}

We found an expression for the solution $y(x)$ of (\ref{eqy}) for $x>0$ and $y>0$. In a similar way, expressions in the other quadrants can be found. However, taking advantage of the discrete symmetries of equation (\ref{eqy}), these solutions can be immediately written down as follows.

Let $\eta\in[1/2,3/2]$, $\al$ with $|\al|<\al_0$.

{\bf{(i)}} Let $y_1=\phi(x;\al,\eta)$ the solution (\ref{formphi}), defined for $x\in[0,\eta]$, with $y_1(\eta)=0$ and $y_1>0$ for $\al>0$. 

Then:

{\bf{(ii)}} the function $y_2=\phi(-x;-\al,\eta)$ is also a solution of (\ref{eqy}), defined for $x\in[-\eta,0]$; we have $y_2(-\eta)=0$ and $y_2\geq 0$ for $\al>0$.

{\bf{(iii)}} The function $y_3=-\phi(-x;\al,\eta)$  is a solution of (\ref{eqy}), defined for $x\in[-\eta,0]$. We have $y_3(-\eta)=0$ and $y_3\leq 0$ for $\al>0$.

{\bf{(iv)}} The function $y_4=-\phi(x;-\al,\eta)$  is a solution of (\ref{eqy}), defined for $x\in[0,\eta]$ and we have $y_4(\eta)=0$ and $y_4\leq 0$ for $\al>0$.

\subsubsection{Matching at the positive $y$-axis}

Let $\eta,\tilde{\eta}\in[ 1/2,3/4]$ and let $y_1(x)=\phi(x;\al,\eta)$ solution of (\ref{eqy}) as in {\bf{(i)}}, for $x\in[0,\eta]$ and $\tilde{y_2}(x)=\phi(-x;-\al,\tilde{\eta})$ solution as in {\bf{(ii)}}, for $x\in[-\tilde{\eta},0]$. 

The following Lemma finds $\tilde{\eta}$ so that $y_1(0)=\tilde{y_2}(0)$, therefore so that $y_1$ is the continuation of $\tilde{y_2}$:

\begin{Lemma}\label{Lyp}

Let $|\al|<\al_0$ with $\al_0$ satisfying (\ref{estimal0}). Let $\eta$ so that $|\eta-1|<c_1|\al|$ with $\al_0$ small enough so that $c_1\al_0\leq 1/2$.

There exists a unique $\tilde{\eta}=\eta+O(\al)$ so that 
\be\label{eqphi1}
\phi(0;\al,\eta)=\phi(0;-\al,\tilde{\eta})
\ee
 and $\tilde{\eta}$ depends analytically on $\eta$ and $\al$ for $|\al|<\al_1$ for $\al_1$ small enough. We have $|\tilde{\eta}-1|\leq c_2|\al |$ for some $c_2>0$ and
\be\label{serteta}
\tilde{\eta}=\eta-{\eta}^{4}c_{{1}}\alpha+2\,{\alpha}^{2}{\eta}^{7}{c_{{1}}}^{2}-
{\alpha}^{3}{\eta}^{10} \left( 4\,c_{{3}}+\frac{21}{2}\,{c_{{1}}}^{3}+10\,c
_{{2}}c_{{1}} \right) +O(\al^4)
\ee
where $c_n=v_n(1)$ with $v_n(\xi)$ given by  (\ref{v12xi}), (\ref{v3xi}).

\end{Lemma}

{\em{Proof.}}

Let
\be\label{defF}
F(\tilde{\eta},\eta,\al)=\phi(0;-\al,\tilde{\eta})-\phi(0;\al,\eta)
\ee
which is a function analytic in $(\tilde{\eta},\eta,\al)$ by Lemma\,\ref{L1} and relation (\ref{defdelta}).

We have $F(\eta,\eta,0)=0$ and 
$$\frac{\partial F}{\partial \tilde{\eta}}(\eta,\eta,0)=-4\eta^3\left[1-v(1;0)\right]=-4\eta^3\ne 0$$
therefore the implicit equation $F(\tilde{\eta},\eta,\al)=0$ determines $\tilde{\eta}=\tilde{\eta}(\eta,\al)$ as an analytic function of $\al$, $\eta$ for $\al$ small. 

We have
$$|\tilde{\eta}-\eta|\leq |\al|\,\sup_{|\al|<\al_1,|\eta|<c_1\al_1}\big|\frac{\partial\tilde{\eta}}{\partial\al}\big|=c_1'|\al| $$
therefore
$$|\tilde{\eta}-1|\leq |\tilde{\eta}-\eta|+|{\eta}-1|\leq (c_1+c_1')|\al |\equiv c_2|\al |$$

The expansion of $\tilde{\eta}$ in power series of $\al$ is found as follows. Using (\ref{formphi}) equation (\ref{eqphi1}) becomes
$$\tilde{\eta}^{4}-\tilde{\eta}^{4}  v
 \left( 1 ;-2\,\tilde{\eta}^{3}\alpha \right)={\eta}^{4}-{\eta}^{4}  v
 \left( 1 ;2\,{\eta}^{3}\alpha \right) $$
 where, substituting (\ref{servxi}) and $\tilde{\eta}=\sum_{n\geq 0}\delta^n\tilde{\eta}_n$ followed by power series expansion in $\al$ we obtain (\ref{serteta}).
\qed

 \subsubsection{Matching at the negative $y$-axis}

Let $\eta,\tilde{\tilde{\eta}}\in[1/2,3/4]$ and $\tilde{\eta}$ given by Lemma\,\ref{Lyp}. Consider the solution $\tilde{y_3}=-\phi(-x;\al,\tilde{\eta})$ as in {\bf{(iii)}}, for $x\in[-\tilde{\eta},0]$, with $\tilde{y_3}(-\tilde{\eta})=0$. Therefore  $\tilde{y_3}$ is the continuation of  $\tilde{y_2}$.

Let $\tilde{\tilde{y_4}}(x)=-\phi(x;-\al,\tilde{\tilde{\eta}})$ be a solution of (\ref{eqy}) as in {\bf{(iv)}}, for $x\in[0,\tilde{\tilde{\eta}}]$.

The following Lemma finds $\tilde{\tilde{\eta}}$ so that $\tilde{y_3}(0)=\tilde{\tilde{y_4}}(0)$, therefore so that $\tilde{\tilde{y_4}}$ is the continuation of $\tilde{y_3}$:

\begin{Lemma}\label{Lyp34}
Let $|\al|<\al_1$ and $|\tilde{\eta}-1|\leq c_2|\al |$ with $\al_1$ small enough so that $\tilde{\eta}\in[1/2,3/2]$.

There exists a unique $\tilde{\tilde{\eta}}>0$ so that $\phi(0;-\al,\tilde{\tilde{\eta}})=\phi(0;\al,\tilde{\eta})$ and $\tilde{\tilde{\eta}}$ depends analytically on $\tilde{\eta}$ and $\al$ for $|\al|<\al_2$ for $\al_2$ small enough. We have $|\tilde{\tilde{\eta}}-1|\leq c_3|\al |$ for some $c_3>0$ and
\be\label{serilteta}
\tilde{\tilde{\eta}}=\tilde{\eta}-\tilde{\eta}^{4}c_{{1}}\alpha+2\,{\alpha}^{2}\tilde{\eta}^{7}{c_{{1}}}^{2}-
{\alpha}^{3}\tilde{\eta}^{10} \left( 4\,c_{{3}}+\frac{21}{2}\,{c_{{1}}}^{3}+10\,c
_{{2}}c_{{1}} \right) +O(\al^4)
\ee
where $c_n=v_n(1)$ with $v_n(\xi)$ given by  (\ref{v12xi}), (\ref{v3xi}).

\end{Lemma}

{\em{Proof.}}

We need to find $\tilde{\tilde{\eta}}=\tilde{\tilde{\eta}}(\tilde{\eta},\al)$ so that 
$$F(\tilde{\tilde{\eta}},\tilde{\eta},\al)=\phi(0;-\al,\tilde{\tilde{\eta}})-\phi(0;\al,\tilde{\eta})=0$$

Note that the function $F$ above is the same as (\ref{defF}). By Lemma\,\ref{Lyp} the present Lemma follows.
\qed

\subsection{The first return map}

Let $\eta\in[1/2,3/2]$ and $\al_2$ as in Lemma\,\ref{Lyp34}. Then $\tilde{\tilde{\eta}}$ given by Lemma\,\ref{Lyp34} is the first return to the positive $x$-axis of the solution with $x(0)=\eta,y(0)=0$ and it is analytic in $\al$ and $\eta$, therefore, by (\ref{defdelta}), it is analytic in $\epsilon$ for fixed $\eta$. 

Combining (\ref{serteta}) and (\ref{serilteta}) we obtain $\tilde{\tilde{\eta}}$ as a convergent power series in $\eta$, with coefficients dependent on $\eta$, whose first terms are
\begin{multline}\label{ret_eta}
\tilde{\tilde{\eta}}=\eta-2\,{\eta}^{4}c_{{1}}\alpha+8\,{\eta}^{7}{c_{{1}}}^{2}{\alpha}^{2
}\\
-{\eta}^{10} \left( 20\,c_{{2}}c_{{1}}+8\,c_{{3}}+49\,{c_{{1}}}^{3}
 \right) {\alpha}^{3}+O \left( {\alpha}^{4} \right) 
 \end{multline}
 
 To obtain the point $\tilde{\tilde{X}}$ where the solution of (\ref{sysXY}) with $X(0)=\epsilon>0$ and $Y(0)=0$ first returns to the positive $X$-axis let $\alpha=2^{1/2}\epsilon^3$ and multiply (\ref{ret_eta}) by $\epsilon$ (since we have $X=\epsilon x$)
and, finally, let $\eta=1$. We obtain that $\tilde{\tilde{X}}$ is analytic in $\epsilon$ for small $\epsilon$ and 
$$\tilde{\tilde{X}}=\epsilon-2^{3/2}\,c_{{1}}{\epsilon}^{4}+16\,{c_{{1}}}^{2}{\epsilon
}^{7}-2^{3/2}\, \left( 20\,c_{{2}}c_{{1}}+8\,c_{{3}}+49\,{c_{{1}}}^{3}
 \right) {\epsilon}^{10}+O \left( {\epsilon}^{13} \right) 
$$
where $c_n=v_n(1)$ with $v_n$ given by (\ref{v12xi}), (\ref{v3xi}).

{\bf{Remark.}} The first coefficient of the return map  (\ref{ret_eta}) is, up to a sign, the Melnikov integral of the system (\ref{syse}), see\cite{Francoise_Iterated_int}; of course, the present results are in agreement with this fact (see \S\ref{Melnikov} for details).

\

\

{\bf{Acknowledgement.}} The author is grateful to Chris Miller for suggesting the problem.

\section{Appendix}\label{Melnikov}

With the notation $H=y^2+x^4$ and $\omega=y^3\, dx$ the Melnikov integral of (\ref{syse}) is the quantity $M(T)=\int_{H=T}\omega$ where $T>0$. If $T$ is the parametrization for the restriction of $H$ to the half-line $x>0$ (which means, in the notations used in the present paper, that $T=\eta^4$) then the return map of (\ref{syse}) has the form $T\mapsto T-\alpha M(T)+O(\alpha^2)$ \cite{Francoise_Iterated_int}.

For the present system we have
\begin{multline}\nonumber
M(T)=\int_{y^2+x^4=T}y^3\, dx=4\int_0^{T^{1/4}}\, (T-x^4)^{3/2}\, dx\\
=4\,T^{7/4}\int_0^1(1-s^4)^{3/2}\, ds=8\, T^{7/4} c_1
\end{multline}
Taking the fourth root in the return map of $T$ we obtain (\ref{ret_eta}).

\end{document}